\magnification 1200

\def\HUP{1}
\def\HUM{2}
\def\COD2{3}
\def\UC{4}
\def\moddeg{5}
\def\GLC{6}

\font\tenmsb=msbm10
\font\sevenmsb=msbm10 at 7pt
\font\fivemsb=msbm10 at 5pt
\newfam\msbfam
\textfont\msbfam=\tenmsb
\scriptfont\msbfam=\sevenmsb
\scriptscriptfont\msbfam=\fivemsb
\def\Bbb#1{{\fam\msbfam\relax#1}}

\def\sqr#1#2{{\vcenter{\vbox{\hrule height.#2pt \hbox{\vrule width.#2pt
height #1pt \kern #1pt \vrule width.#2pt}\hrule height.#2pt}}}}
\def\qed{$\  \sqr74$}

\def\mapdown#1{\Big\downarrow \rlap{$\vcenter{\hbox{$\scriptstyle#1$}}$}}
\def\mapup#1{\Big\uparrow \rlap{$\vcenter{\hbox{$\scriptstyle#1$}}$}}
\def\lmapdown#1{\llap{$\vcenter{\hbox{$\scriptstyle#1$}}\Big\downarrow $}}
\def\mapright#1{\smash{\mathop{\longrightarrow}\limits^{#1}}}

 \def\la{\longrightarrow}

 \def\ni{\noindent}
 \def\cl{\centerline}
\def\df{\noindent {\bf Definition.\  }}

 \def\pf{\noindent {\it Proof.\  }}

 \def\s{\sigma}
\def\a{\alpha}

\def\g{\gamma}

\font\gothicf=eufm10
\font\sgothic=eufm7
\font\ssgothic=eufm5
\textfont5=\gothicf
\scriptfont5=\sgothic
\scriptscriptfont5=\ssgothic

\def\C{{\Bbb C}}

\def\P{{\Bbb P}}
\def\Q{{\Bbb Q}}

\def\Z{{\Bbb Z}}
\def\N{{\Bbb N}}

\def\Im{{\rm Im}}
\def\codim{{\rm codim}}


\def\NL{N(L,g)}

\def\Mg{{\overline M_g}}

\def\FBS{F_g(B,S)}
\def\FVT{F_g(V,T)}


\cl{{\bf REMARKS ABOUT UNIFORM BOUNDEDNESS}}
\cl{{\bf  OF RATIONAL POINTS OVER FUNCTION FIELDS}}

\

\cl{LUCIA CAPORASO}

\
\

\ni
{\bf 1. Introduction and preliminaries.} A curve $X$ of genus at least $2$ defined
over a function field $L$ has only finitely many $L$ rational ponts, unless
it is isotrivial. Similarly,   a curve of genus at least $2$ defined over a number
field $F$ has a finite set of
$F$-rational points. These well known facts are celebrated theorems of Y.Manin and
G.Faltings, originally conjectured by L.J.Mordell and S.Lang.

We study here questions of uniformity for the cardinality of such sets of rational points, in the
function field case.
For number fields, there are a number  of open conjectures, such as the following
(Uniform Mordell Conjecture for number fields):
{\it {Fix $g\geq 2$ and a number field $F$; there exists a number $B_g(F)$ such that  any curve of genus
$g$ defined over $F$ has at most $B_g(F) $ rational points  over $F$.}}
Interest in such problems was revived after it was proved in [CHM] that the 
Conjecture above is a consequence of a famous, open, Conjecture (usually attributed to
S.Lang and E.Bombieri) on the non-density of rational points in varieties of general
type (see  also [A], [AV] and [Pa]). 

In this paper we investigate similar issues for curves over function fields.
Some partial results were obtained in [Mi] and in [C] 
were the existence of uniform bounds for the sets of rational points
is established.
Such bounds depend on suitable numerical invariants of the function field, on
the
genus $g$ of the curves and on the
degree of the  locus of bad reduction (that is, the locus of singular fibers).

We shall also study here the strictly related ``uniform Shafarevich problem";
a famous Theorem of 
A.N.Parshin and 
S.Ju.Arakelov ([A] and [P]) states that
{\it if $B$ is a smooth complex curve and $S\subset B$ a finite subset, then there
exists only a finite number of non-isotrivial families of smooth curves of fixed genus
$g\geq 2$ over $B-S$}.
Parshin first proved it under the assumption that $S=\emptyset$; 
Arakelov generalized it a few years
later. In [P]  Parshin shows also that the Theorem above  implies finiteness of
rational points for non-isotrivial curves of genus at least $2$,
providing the above mentioned link between the Shafarevich problem and the Mordell
problem.  Recall  that his argument,  known as the ``Parshin trick", is
valid  for both number  fields and  function fields.

A first uniform version of the 
Theorem of Parshin and Arakelov above is obtained in
[C]. We  here generalize it by a stronger uniform result
valid for families of curves over bases of any dimension.
This is done in Section 2, where we obtain bounds 
(for the sets of curves with fixed degeneracy locus as well as for the sets of rational
points) that only depend on the degree of a polarization on the base variety, and on
the degree of the locus of bad reduction.
A  stronger  result can be obtained for curves having good reduction in
codimension 1 (Theorem \COD2).
In Section 3 we will consider families with maximal variation of moduli, 
using the geometry of the moduli space of curves to approach our problems.

\

  We work over
$\C$, by $V$ we shall denote a smooth , irreducible, projective variety over $\C$,
whose field of rational functions will be
$L := \C (V)$. 
Special interest will be given to varieties of dimension $1$, for which we shall use
the following notation: $B$ is a smooth irreducible curve and $K$ its field of
rational functions. We fix integers $q\geq 0$, $g\geq 2$ and $s \geq 0$ throughout.
The genus of $B$ will be denoted by $q$.

We shall consider smooth curves of genus $g$ over the function field $L$ (or $K$), 
which
 can also be viewed as  families of  curves  over  $V$, 
such that there is a non-empty  open subset
of
$V$ over which  the fibers are all smooth. We shall  always
assume that such a family (or curve) is not isotrivial, i.e. the smooth fibers are not all isomorphic.

To be more precise, we introduce the following sets:
let $B$ be a fixed curve and let $S\subset B$ be a  finite set of points.

\

\ni
\df
 $\FBS$ shall denote the set of equivalence classes of non-isotrivial families
$f:X\la B$ such that $X$ is a smooth relatively minimal surface and the fiber $X_b$
over every
$b\not\in S$ is a smooth curve of genus $g$. Two such families $f_i:X_i\la B$ for
$i=1,2$ are  equivalent if there is a commutative diagram 
$$
\matrix{X_1&\mapright {\a '} & X_2\cr
\  \  \lmapdown{f_1}&&\mapdown{f_2} \cr
B & \mapright{\a }& B}
$$
where the two horizontal arrows are birational maps.

\

\ni
Using a different terminology, $\FBS$ is the set of $K$-isomorphism classes of non-isotrivial curves of
genus $g$ over $K$, having good reduction outside of $S$.
The Theorem of 
Parshin and 
Arakelov says that
 $\FBS$ is finite. Theorem 3.1 of [C] states that there exists a number $P(g,q,s)$ such
that
$|\FBS |\leq P(g,q,s)$ for every curve $B$ of genus $q$ and for every subset $S$
having at most $s$ points.
We show here (in the end of Section 2) that this result is sharp in the sense that such
a bound must depend on
$s$.

We are interested in function fields of higher transcendence degree.
We can generalize the definition of $\FBS$ as follows. Let
$T\subset V$ be a closed subscheme.

\

\df
 $\FVT$ shall be the set of equivalence classes of non-isotrivial families of smooth
curves of genus
$g$ over $V - T$ 
(the equivalence relation is the same as above, with $B$ replaced by $V-T$).

\

\ni
By the existence and unicity of minimal models for smooth surfaces,  
this definition coincides with the
previous one if $\dim V=1$.
It follows from the results in [C] (3.4) that $\FVT$ is finite.
Our best result on $\FVT$ is Theorem \HUP.

If $X$ is a curve defined over a field $L$, we shall denote by $X(L)$ the set of its
$L$-rational points. If $X$ has genus at least $2$ and it is not isotrivial,
the Theorem of Manin says that $X(L)$ is finite.
Consider now the Uniformity Conjecture for rational points over function fields,
which can be stated  as its arithmetic analogue:
{\it Let $L$ be a function field over $\C$ and let $g\geq 2$ be an integer.
There exists a number $N_g(L)$ such that for every non-isotrivial curve $X$ of genus
$g$ defined over $L$ we have  $|X(L)\leq N_g(L)$.}
For results relating it to the Lang Conjectures about the
distributions of rational
points on varieties of general type, see the work of D. Abramovich and J.F. Voloch
[AV].

Such a Conjecture remains open; our results in that direction are 
Theorems \HUM  \  and  \COD2 \ and Proposition \UC. 

A final piece of notation. $M_g$ denotes the moduli variety of smooth curves of genus
$g$ and $\Mg$ its compactification via Deligne-Mumford stable curves. They are both
integral, normal varieties of dimension $3g-3$.
A universal curve exists only on a proper open subset of $M_g$ (and of $\Mg$).
In particular, a morphism $\phi :Z\la M_g$ does not necessarily come from a family of
curves over $Z$. If this is the case, that is, if there exists a family of smooth
curves $X\la Z$ such that for every $z\in Z$, $\phi (z)$ is the isomorphism class of
the fiber of
$X$ over $z$, we shall say that $\phi$ is a {\it moduli map}.

\

\ni
{\bf   2. Uniformity results for function fields of high transcendence degree.}
We start by a uniform generalization fo the Theorem of Parshin and Arakelov. The
result below is  a strenghtening of 3.4 and 3.5 in [C], in fact the bound $H$ here is
independent of the dimension of $V$ and of
$r$. Such an improvement is obtained by a small technical modification of the
methods in [C].

Notice that statement below remains true if $V$ is replaced by an integral, possibly
singular, projective variety. The proof is essentially the same.

\proclaim Theorem \HUP.
Let $g\geq 2$, $d\geq 1$ $s\geq 0$ be fixed integers. 
There exists a number $H(g,d,s)$ such that for any
smooth, irreducible variety $V\subset \P ^r$ of degree $d$, 
for any closed subscheme $T\subset V$ of degree
$s$, we have $|F_g(V,T)|\leq H(g,d,s)$.
Moreover, if $T$ has codimension at least $2$ in $V$, then the bound $H$ does not depend on $s$.

\pf
{\it Step  1. Slicing $V$ into curves  of bounded genus.}
 Considering
one-dimensional hyperplane sections of
$V$,  we see that  $V$ can be covered by 
smooth curves of degree
$d$ passing through any of its points;  it is a well known fact that the
 genus of a   curve of degree $d$ in projective
space is at most equal to
${d-1 \choose 2}$:
just project the curve birationally onto  a curve of degree $d$ in $\P ^2$. Let 
$$q=q(d) ={d-1 \choose 2}$$ 
so
that
$V$ is covered by curves of geometric genus at most $q$.

{\it Step\  2.  Uniform boundedness of moduli maps.} By Theorem 3.1 in [C], for any
fixed $g, q', s'$ there exists a number 
$P(g,q',s')$ such that for any smooth curve $B$ of  genus $q'$, for any
subset $S\subset B$ of at most $s'$ points, we have that
$|\FBS | \leq P(g,q',s')$.

Define
$$H'=\max_{q'\leq q, s'\leq s}P(g,q',s')$$
 so that $H'$ only depends on $g,d,s$;  let
$U=V-T$.

We claim that $U$ has at most $H'$ moduli
maps to
$M_g$, that is, we claim that there exist at most
$H'$ non-constant, (regular) morphisms $\phi :U\la M_g$ such that there exists a 
(not necessarily unique, see below) family
of smooth curves over $U$ whose moduli map is $\phi$. 
By contradiction, let $n> H'$ and let us assume that there exist $\phi _1,....,\phi _n$ distinct such
moduli maps $\phi _i:U\la M_g$. Let $X_i\la U$ be a non-isotrivial 
family of smooth curves corresponding to $\phi _i$
(since $\phi_i$ is a moduli map, such a family exists, but it is not
necessarily unique). Let
$U'\subset U$ be the non-empty open subset where $\phi _i(u)\neq
\phi _j (u)$ for every $u\in U'$ and for every pair of distinct $i$, $j$. 
Let $p\in U'$ and let $F_i = \phi _i ^{-1}\phi _i (p)$; since $\phi _i$
is not constant, its fiber $F_i$ through $p$ is a proper closed subset of $U'$, 
therefore there exists a
curve
$B\subset V$ of genus at $q'\leq q$ such that $p\in B$ and such that $B\not\subset F_i$
for every $i=1,....,n$; thus the restriction of
$X_i$ to
$B$ is not isotrivial for every
$i$. Let $S=(B\cap T)_{red}$.
Let $Y_i\la B$ be the smooth relatively mininimal completion over $B$ of the restriction of $X_i$ to
$B$. By construction, $Y_1,....,Y_n$, are different elements of $\FBS$, 
which is a contradiction, since
$\FBS$ has at most $P(g,q',s') \leq H' <n$ elements.
This proves the claim. Notice that if $T$ has codimension at least $2$ in $V$ we can 
always choose our $B$ 
so that it does not intersect $T$ at all, and hence $S$ can be taken to be the empty set and $H'$ does
not depend on $s$.

{\it Conclusion.} Given a moduli map $\phi :U\la M_g$ the set of
families that have $\phi$ as moduli map is uniformly bounded, in fact it is bounded
above by a function of $g$ only (see [C] Lemma 3.3), hence we are done.
\qed

\

A  similar argument yields the following uniformity statement for rational points,
stronger than 4.3 and 4.4 in [C]:

\proclaim Theorem \HUM.
Let $g\geq 2$, $d\geq 1$ $s\geq 0$ be fixed integers. There exists a number $N(g,d,s)$ such that for any
smooth, irreducible variety $V\subset \P ^r$ of degree $d$, for any closed subscheme
$T\subset V$ of degree
$s$ and for any non-isotrivial curve $X$ of genus $g$ defined over $L=\C (V)$ and having good reduction
outside of $T$, we have $|X(L)|\leq N(g,d,s)$.
Moreover, if $T$ has codimension at least $2$ in $V$, then the bound $N$ does not depend on $s$.

\pf
{\it{Step  1.}}   Repeat word by word Step 1 in the proof of the
previous Theorem.

{\it{Step  2.}}   Theorem 4.2 in [C]  says that if $g, q', s'$ are fixed
non-negative integeres, there exists a number $M(g,q',s')$ such that for any curve $B$
of genus $q'$, for any subset $S$ of at most $s'$ points in $B$, for any curve $X_B \in
F_g(B,S)$ we have that
$$
|X_B(\C (B))|\leq M(g,q',s').
$$
Arguing as in the proof of 4.4 of [C] one gets that defining 
$$
N(g,d,s):=\max_{q'\leq q, s'\leq s}M(g,q',s').
$$
will suffice to our statemnt.
\qed

\

To conclude, we show that for curves having good reduction in codimension 1, stronger
finiteness results hold.
Let $L$ be a function field over $\C$ and let $V$ be a smooth, projective, complex
variety of positive dimension such that
$L = \C (V)$.

\

\df
Let $C_g^2(L) $  be  the set of $L$-isomorphism classes 
of non-isotrivial curves of
genus
$g$ over
$L$ having good reduction in codimension 1.

\

\ni
In other words, $C_g^2(L)$ is the set of equivalence classes of non-isotrivial
families $X\la V$ of curves of genus $g$ over $V$ such that there exists a closed
subscheme
$T\subset V$ of codimension at least $2$ with the property that $X_v$ is smooth for
every $v \not\in T$.

\proclaim Theorem \COD2.
\item  a\negthinspace)  $C_g^2(L)$ is finite 
\item  b\negthinspace)  There exists a number $N_g^2(L)$ such that for every curve
$X\in C_g^2(L)$  we have
$|X(L)|\leq N_g^2(L)$.

\pf
We shall use moduli maps. 
Denote by $
 M_g^2(L)
$
the set of equivalence classes of non-constant rational maps $\phi : V \la M_g$ such
that there exists an open subset
$U^{\phi }\subset V$ with the following properties:

\ni
1. The complement of $U^{\phi }$  has codimension at least $2$ in $V$,

\ni
2. $\phi$ is regular on $U ^{\phi }$,

\ni
3. There exists a (non-isotrivial) family of smooth curves of genus $g$ over $U
^{\phi}$ such that
$\phi$ is its moduli map.

\ni
4. Two such maps $\phi$ and $\psi$ are equivalent iff they coincide on some (non-empty)
open subset of $V$.

\

There is a natural surjective map of sets:
$$
\mu : C_g^2(L) \la M^2_g(L)
$$
sending a curve over $L$ to its moduli map (it is easy to see that $\mu$ is well
defined). Now, $\mu$ has finite fibers (Lemma 3.3 in [C]) and is surjective by
definition. Thus $C^2_g(L)$ is finite if and only if $M^2_g(L)$  is finite.

Part b) is an immediate consequence of part a), by the Theorem of Manin.
We will prove our result by showing that $M^2_g(L)$ is finite by induction on $\dim V$.
If $\dim V =1$, the the finiteness of $C^2_g(L)$ and  of $M^2_g(L)$ is the Theorem of
Parshin (the locus of bad reduction being empty in such a case). Let then $\dim V \geq
2$ and suppose that
$M^2_g(L)$ is infinite. Notice that $M^2_g(L)$ is dominated by a union
of finite sets as follows; if $T$ is a closed subset of $V$, 
denote by $M_g(V,T)$  the set of equivalence classes of moduli maps to $M_g$
that are regular on  $V-
T$, then $M_g(V,T)$ is finite, by Theorem \HUP 
  \ and  Lemma 3.3 in [C].   We have  a
natural, surjective map
$$\bigcup _{\codim _V T\geq 2} M_g(V,T) \la M^2_g(L)
$$
hence, if $M^2_g(L)$ is infinite, so is  the union on the left hand side.  Then
there exists a countable collection $\{ T _n,\  n\in \Z\}$, with
$T_n$  a closed subset of $V$ of codimension at least $2$, such that the set
$$
M := \bigcup _{n\in \Z} M_g(V,T_n)
$$
is infinite.
Now, $M$ itself being a countable set,  we shall put an ordering on it:
$$M = \{ \phi
^i,\  i \in \N \}$$
For every pair of distinct $i,j$, denote by $U^{i,j}$ the non-empty open subset
of $V$ such that $U^{i,j} \subset U^{\phi ^i}\cap U ^{\phi ^j}$ and 
$\phi ^i(u) \neq \phi ^j(u)$ for every $u\in U^{i,j}$.
The $U^{i,j}$s form a countable collection of non-empty open subsets of $V$, whose
intersection $I$ is dense in $V$. Let $p\in I$ and, for every $i \in \N$,  let $F_i =
(\phi ^i)^{-1}\phi ^i(p)$  be the fiber of $\phi ^i$ through $p$. Since $\phi ^i$ is
non-constant (by assumption)
$F_i$ is a proper  closed subset of $V$, thus the complement of $\cup _{i \in \N}F_i$
intersects $I$ in a subset $J$, with $J$ dense in $V$.
Fix a (non-degenerate) projective model of $V$ in some projective space. Then there
exists a hyperplane
$H$ such that $p\in H$, such that $H\cap J  \neq \emptyset$ and such that $H$ does not
contain any $T_n$. Let $W=H\cap V$, we can furthermore choose $H$ so that $W$ is
smooth.
By construction we have 

(a) $\dim W = \dim V -1$

(b) $\dim T_n \cap W = \dim T_n -1 \leq \dim W -2$

(since $H$ does not contain any $T_n$)

(c) $\forall \phi ^i \in M$, the restriction $\phi ^i_{|W}$ is not constant

(since $p\in W$ and $W\cap J \neq \emptyset$)

(d) $\forall i\neq j$ we have $ \phi ^i _{|W}\neq \phi ^j_{|W}$

(since $W\cap I \neq \emptyset$)

\ni
hence $\phi ^i_{|W}\in M^2_g(\C (W))$ and the restriction to $W$ gives  
an inclusion (by (d) above)
$M\hookrightarrow M^2_g(\C (W))$.
Thus $M^2_g(\C (W))$ is infinite. This is a contradiction with the inductive
assumption.
\qed 

\

See [Md]
for  an analogue over $\Q$.
Part a) of this result should be compared with the examples of A. Beauville (in [B],
section 5) or with the example  below. They show that the
assumption that the curves have good reduction in codimension 1 is crucial,
that is, a) is false without that assumption.
On a different vein, compare also with Proposition 6.
The example that we are going to describe shows that there is no hope of getting a
substancially stronger
uniform version of the Shafarevich Conjecture for function fields; in other words,
any uniform bound on $|\FBS |$ must depend on the degree of $S$. 

What happens to the cardinality of $\FBS$ when $s$ grows
while $g$ and $q$ (or even $B$) stay fixed? The way we defined $\FBS$, 
it is an exercise to show that
its cardinality is not bounded; but this is just because the families
parametrized by
$\FBS$ are not  required to
have singular fibers over $S$.
The  interesting question is about the asymptotics of the cardinality of that subset
of $\FBS$ parametrizing families of curves that have singular fiber over every point
of $S$. We will make this precise now,  describing an example 
suggested by  J.DeJong, showing  that the set of fibrations with fixed degeneracy
locus is not bounded, as the cardinality of the degeneracy locus grows. 

Fix  $g\geq 2$ and  $B=\P ^1$, given  a subset $S\subset \P ^1$ denote by $F(S)\subset F_g(\P ^1, S)$ the
set of all genus
$g$ non-isotrivial fibrations $X\la \P ^1$    such that  the fiber 
$X_b$ is smooth if and only if $b\notin S$.

Let $S=\{a_1,....a_s\}$ be a set of generic points in $\P ^1$, and let $I\cup J =\{ 1,2,....,s\}$ be a
partition of
$\{ 1,2,....,s\}$ in two disjoint subsets such that $|I|=5$.
Define a non-isotrivial fibration $X_I$ of curves of genus $2$ over $\P ^1$ by the affine equation
$$
y^2=(x-t)\Pi _{i\in I}(x-a_i)\Pi _{j\in J}(t-a_j)
$$
with $t$ affine coordinate in $\P ^1$. For $t\notin S$ (and $t \neq \infty$ ) we get a smooth curve of
genus
$2$. For
$t=a_i$ with $i\in I$,  we get a nodal curve and for $t=a_j$, $j\in J$, we get a singular, non-reduced
curve. Thus $X_I\in F(S\cup \infty )$ and by varying the partition $I\cup J$ we get a total of $s \choose
5$ different such fibrations.
Hence the cardinality of $F(S)$ goes to infinity, as $|S|$ grows.

One final word about this example. 

First, we make two comments:
the given family has fibers of genus $2$, but of course   one can construct the same example for any
genus, (just replace the integer 5 by a higher odd number), obtaining  families of hyperelliptic
curves. 

The second comment is about the
 singular fibers over $a_j$ with $j\in J$, which are  not stable curves;  
their semistable reduction is actually a
smooth curve. The remaining $5$ fibers over $a_i$ are instead nodal.
In other words,  the  moduli map
$\phi _I$ associated to the family
$X_I\la
\P ^1$
$$
\phi _I : \P ^1 \la \overline{M_2}
$$
(such that $\phi _I(t)$ is the isomorphism classes of the fiber of $X_I$ over $t$) intersects the
boundary $\Delta _2$ in exactly 5  points, regardless of the cardinality of $S$.

We ask:

(a) Can one find  similar examples whose fibers  do not belong to
any proper closed subset of $M_g$? 

(b) Is the same ``unboundedness" result true for families of stable curves? In other words,
does there exist a similar  example all of whose singular fibers are nodal?

\

\ni
{\bf 3. Uniformity for ``truly varying" curves.}
This section contains results that are independent of the degeneracy locus.
Given a family $X\la V$ of generically smooth curves of genus $g$ over $V$,
we get a natural rational map $\phi:V\la M_g$ (regular on a non-empty open subset  of $V$).
The dimension of the image of $\phi$ is called the {\it variation of moduli} of the family;
we shall say that the family has {\it maximal variation of moduli} if 
$$
\dim \Im \phi =\min \{\dim V, 3g-3\}.
$$
We shall say that a curve over $L=\C (V)$ has maximal variation of moduli if
a corresponding family of curves over $V$ does.

Thus the condition of having maximal variation of moduli can be 
interpreted as saying that the family (or the curve) is {\it truly varying} and can be
viewed as a generalization of the
non-isotriviality  condition. Obviously, 
if the base field has transcendence-degree $1$,  a curve is
non-isotrivial if and only if it has maximal variation of moduli.

\

\df
Let $L$ be a function field, we define $C_g(L)$ to be the set of 
$L$-isomorphism classes of curves of genus $g$ defined over $L$
and having maximal variation of moduli.

\

\proclaim Proposition \UC. Let $g\geq 24$ and let $L$ be a  function field  
of transcendence degree $3g-3$.
Then 
\item   a \negthinspace) $C_g(L)$ is finite 
\item  b\negthinspace)  There exists a number $\NL$ such that for every curve $X$ of
genus
$g$ defined over $L$ and having maximal variation of moduli, we have $|X(L)|\leq \NL$
\item c\negthinspace)  There exists a function
$P_g(n, m)$ such that for every  $V$ of general type,
we have $|C_g(L)|\leq P_g(\dim V, K_V^{\dim V})$.

\pf
The assumption $g\geq 24$ implies that $\Mg$ is of general
type (for this famous result of J.Harris and D.Mumford  we refer to [HMu] and to
 6F in [HM]).  

Denote by $R(V,\Mg )$ the set of dominant, rational maps from $V$ to $\Mg$. A Theorem
of Kobayashi-Ochiai [KO] implies that, $\Mg $ being of general type,
$R(V,\Mg )$ is finite. 
Notice now that there is a natural bijection between $C_g(L)$ and $R(V,\Mg )$:
to a truly varying curve $X$ of genus $g$ over $L$ we can associate its moduli map
$\phi _X\in R(V,\Mg )$. The fact that such a correspondence is bijective follows from
the existence of the universal curve over an open subset of $\Mg$. Thus 
$C_g(L)$ is finite. 

By the Theorem   of Manin, any curve in $C_g(L)$ has a finite set of $L$-rational
points, thus part b)  follows immediately from  a)

Part c) is proved like part a); we can in this case apply a strengthening of the
Theorem of Kobayashi-Ochiai provided by T. Bandman and D. Markushevich. From [BM] we
obtain that, $V$ and $\Mg$ being of general type and $\Mg$ having canonical
singularities (Theorem 1 in [HMu]), there exists a function of $g$,
of $\dim V$  and of
$K_V^{\dim V}$ bounding the cardinality of $R(V,\Mg )$ and hence that of $C_g(L)$.
\qed

\

Let $u:{\cal C}_g\la M_g^o$ be the universal curve over the moduli space of automorphism free smooth
curves of genus $g$, so that the fiber of $u$ over the point corresponding to the
curve $X$ is 
$X$ itself. 

It is a well known fact  (see [HM] 2D) that $u$ has no rational sections, thus, ${\cal
C}_g$ has no rational point over the function field of $M_g$.
In fact much more is known: the Picard group of ${\cal C}_g$  is generated over the
Picard group of $M_g$ by the relative dualizing sheaf $\omega _u$; therefore a
multisection of $u$ must have degree over $M_g^o$ equal to a multiple of $2g-2$. 

We apply this to obtain that if $V$ is a variety of dimension $3g-3$ and $X$ is a
curve of genus $g$ over $L$ having maximal variation of moduli, then a necessary condition for $X$
to have a rational point over $L$ is that its  moduli map have degree equal to a multiple of $2g-2$.
 This follows easily by looking at the  commutative diagram
$$
\matrix{X\  &\mapright{\gamma }&{\cal C}_g\cr
\  \lmapdown{f }  \mapup{\s }& &\mapdown {u}\cr
V \  &\mapright {\phi _X} &M_g }
$$
where the horizontal arrows are rational maps and 
$\s$ is the rational  section corresponding to a rational point of $X$ over $L$. Let $\tau = \g \circ
\s:V\la {\cal C}_g$ and let
$\rho : \Im \tau \la M_g$; by what we said, $\deg \rho = n(2g-2)$ for some integer $n$.
We finally obtain
$$
\deg \phi _X = \deg \tau \cdot \deg \rho = m(2g-2).
$$
Where by $\deg \phi _X$ we mean the degree of the restriction of $\phi _X$ to the non-empty
open subset of $V$ where $\phi _X$ is a regular and finite map.
Let us call such a number $\deg \phi _X$   the {\it modular degree} of a family
$X\la V$; this definition is general, provided that $X\la V$ has maximal variation
of moduli and that $\dim V \leq 3g-3$.
We just proved the following

\proclaim Lemma \moddeg. Let $V$ be a variety of dimension $3g-3$ with function field
$L$ and let $X$ be a smooth curve of genus $g$ over $L$ having maximal variation of
moduli. Then either
$X(L)=\emptyset$ or the modular degree of $X$ is a multiple of $2g-2$.

 \

\ni
The following
well known conjecture is open:

\proclaim Geometric Lang Conjecture.
Let $W$ be a variety of general type defined over $\C$. Then there exists a proper
closed subvariety $Z_W$ of $W$ containing all positive dimensional subvarieties of $W$
that are not of
general type.

In particular, according to such a Conjecture, 
all curves in $W$ having genus at most $1$ are contained in $Z_W$.

Consider now $\Mg$, and let $Z_g\subset \Mg$  be defined as the closure of the union of
all integral curves in $\Mg$ having geometric genus at most equal to $1$.
Since $\Mg $ is of general type if $g\geq 24$, the above
conjecture would imply that $Z_g$ is a proper, closed subset of $\Mg$ for all $g\geq 24$.

As a consequence, we get the following

\proclaim Lemma \GLC. Let $g\geq 24$ and $B$ be a curve of
genus $q$.
The Geometric Lang Conjecture implies that if $X\la B$ is a non-isotrivial family of curves of
genus $g$, passing through the general point of $M_g$, then the modular degree of $X$ is at most
$q-1$.

\pf
As we mentioned above, the union of all curves in $\Mg $ of genus at most $1$  is contained in a
proper closed subset $Z_g$  of $\Mg$. The condition that the given family of curves goes through
the general point of $\Mg$, combined with the Geometric Lang Conjecture, implies that $\Im \phi _X
\not\subset Z_g$. Thus the geometric genus of $Im \phi _X$ is at least $2$.
By the Riemann-Hurwitz formula, the degree $d$ of  a dominant map of a curve $B$ of genus $q$
onto a curve $C$ of geometric genus  $p\geq 2$ is at most equal to $q-1$, in fact the formula
gives
$$
d= {2q-2 - r \over 2p-2} \leq {q-1 \over p-1}\leq q-1
$$
since $r\geq 0$ (being the degree of the ramification divisor) and $p\geq 2$ by assumption.
\qed

\

\

\ni
{\bf Aknowledgments} 
I am grateful to Olivier Debarre and to Johann DeJong for useful conversations and
to Felipe Voloch for indicating  relevant references.
Special thanks to Dan Abramovich for pointing out a  
serious mistake in a previous version of
this paper. 

\

\ni
{\bf References}

\

\ni
[A] D.Abramovich. {\it 
On the number of stably integral
points on an elliptic curve.} Inventiones Math. 128, 481-494 (1997).

\ni
[AV] D.Abramovich, J.F.Voloch:  {\it Lang's conjectures, fibered powers, and
uniformity.} New York J. Math. 2 20-34 (1996)

\ni
[A] S.Ju.Arakelov: {\it Families of algebraic curves with fixed degeneracies.} Izv.
Akad. Nauk. SSSR Ser. Mat.35 (1971) no. 6, 1277-1302.

\ni
[BD] T.Bandman,  G.Dethloff:
{\it Estimates of the number of rational mappings from a fixed variety to varieties of general
type.}  
Ann. Inst. Fourier (Grenoble) 47 (1997), no. 3, 801-824.  

\ni
[BM] T.Bandman,  D.Markushevich: 
{\it On the number of rational maps between varieties of general type.}   J. Math.
Sci. Univ. Tokyo 1 (1994), no. 2, 423-433. 

\ni 
[B] A.Beauville: {\it Expos\'e n. 6} in  {\it  S\'eminaire sur les pinceaux des
courbes de genre au moins deux.} Ast\'erisque 86 (1981).

\ni
[C] L.Caporaso: {\it  On certain uniformity properties of curves over function fields.} 
Preprint, AG/9906156, to appear in Compositio Mathematica.

\ni
[CHM] L.Caporaso, J.Harris, B.Mazur: {\it Uniformity of rational points.} 
Journal of American Mathematics
Society, Vol. 7, N. 3, January 97 
p. 1-33.

\ni
[DF] M.De Franchis: {\it Un teorema sulle involuzioni irrazionali.} Rend.Circ.Mat Palermo 36 (1913), 368.

\ni [HM] J.Harris, I.Morrison: {\it Moduli of curves.} Graduate texts in Math. 187 (1998) Springer

 \ni  [HMu] J.Harris,  D.Mumford: {\it On the Kodaira dimension of the
moduli space of curves.}  Inventiones Math. 67
(1982), no. 1, 23--88.

\ni
[KO] S.Kobayashi,\  T.Ochiai: {\it Meromorphic mappings into compact complex spaces of general type.}
Inventiones Math. 31
pp. 7-16 (1975).

\ni
[Md] M.Martin-Deschamps: {\it Conjecture de Shafarevich pour les corps de fonctions
sur
${Q}$.} In  Ast\'erisque No. 127, Appendice \`a l'expos\'e IX
(1985), 256-259. 

\ni
[Ma] Y.Manin: {\it 
Rational points of algebraic curves over function fields.} Izv. Akad. Nauk. 27 (1963),
1395-1440. 

\ni
[Mi] Y.Miyaoka: {\it Themes and variations on inequalities of Chern classes.}
Sûgaku 41 (1989), no. 3, 193-207. 

\ni
[Pa] P.Pacelli. {\it Uniform boundedness of rational points.} Duke Math.J. Vol 88, No
1 pp 77-102 (1997). 

\ni
[P] A.N.Parshin: {\it Algebraic curves over function fields.} Izv. Akad. Nauk. SSSR Ser. Mat.32 (1968) no. 5,
1145-1170.

\ni
[S] L.Szpiro: {\it Expos\'e n. 3} in  {\it S\'eminaire sur les pinceaux des courbes de genre au moins deux.}
Ast\'erisque 86 (1981).

\ni
[V] E.Viehweg: {\it Quasi-projective moduli for polarized manifolds.} Ergebnisse der
Mathematik (30) Springer.

\

\ni
Address of Author:

\ni
L.Caporaso: Universit\`a degli Studi del Sannio, Benevento, Italy

and Massachusetts Institute of Technology, Cambridge, MA, USA

caporaso@math.mit.edu

\

\

\ni
(Revised December 2000)

\end

\end